\documentclass[11pt]{article}

\usepackage{mathrsfs}
\usepackage{amssymb}
\usepackage{amscd}
\usepackage{amsfonts}

\font\tenmsb=msbm10
\font\sevenmsb=msbm7
\font\fivemsb=msbm5

\catcode`\@=11
\ifx\amstexloaded@\relax\catcode`\@=\active
\endinput\else\let\amstexloaded@\relax\fi
\def\spaces@{\space\space\space\space\space}
\def\spaces@@{\spaces@\spaces@\spaces@\spaces@\spaces@}
\def\space@.{\futurelet\space@\relax}
\space@.
\def\Err@#1{\errhelp\defaulthelp@\errmessage{AmS-TeX error: #1}}
\def\relaxnext@{\let\next\relax}
\def\accentfam@{7}
\def\noaccents@{\def\accentfam@{0}}
\def\Cal{\relaxnext@\ifmmode\let\next\Cal@\else
\def\next{\Err@{Use \string\Cal\space only in math mode}}\fi\next}
\def\Cal@#1{{\Cal@@{#1}}}
\def\Cal@@#1{\noaccents@\fam\tw@#1}
\def\Bbb{\relaxnext@\ifmmode\let\next\Bbb@\else
\def\next{\Err@{Use \string\Bbb\space only in math mode}}\fi\next}
\def\Bbb@#1{{\Bbb@@{#1}}}
\def\Bbb@@#1{\noaccents@\fam\msbfam#1}
\newfam\msbfam
\textfont\msbfam=\tenmsb
\scriptfont\msbfam=\sevenmsb
\scriptscriptfont\msbfam=\fivemsb

\def\NN{{\Bbb N}}
\def\ZZ{{\Bbb Z}}
\def\RR{{\Bbb R}}
\def\QQ{{\Bbb Q}}

\def\SS{{\Bbb S}}
\def\o{{\omega}}
\def\O{{\Omega}}

\def\t{{\theta}}

\def\l{{\lambda}}

\def\pa{{\partial}}

\def\beq{\begin{equation}}
\def\eeq{\end{equation}}

\def\qedbox{$\rlap{$\sqcap$}\sqcup$}
\def\skipaline{\removelastskip\vskip12pt plus 1pt minus 1pt}

\def\Proof{\removelastskip\skipaline
\noindent \it Proof. \rm}

\newtheorem{Theorem}{Theorem}
\newtheorem{Lemma}{Lemma}[section]

\newtheorem{Corollary}{Corollary}[section]
\newtheorem{Remark}{Remark}[section]

\@addtoreset{equation}{section}

\begin{document}

\title{Boundedness for Second Order Differential Equations with Jumping p-Laplacian and an oscillating term
\author{Xiao Ma \footnote{  Email: mx5657@sina.cn}
, \ \ Daxiong  \ Piao \footnote{The corresponding author.Supported by the NSF of Shangdong Province(No.ZR2012AM018), Email: dxpiao@ouc.edu.cn}
\\
        School of Mathematical Sciences,\ Ocean University of China,\\
        Qingdao 266100,\ \ P.R. China
\\
        Yiqian Wang \footnote{Supported by the NNSF of Chian(No.11271183),Email:yqwangnju@yahoo.com}
\\
      Department of Mathematics,\ Nanjing University
\\
        \ \ Nanjing 210093,\ \ P.R. China
}
\date{}}

\maketitle

\begin{abstract}In this paper, we are concerned with the boundedness of all the solutions for a kind of second order differential equations with p-Laplacian and an oscillating term $(\phi_p(x'))'+a\phi_p(x^+)-b\phi_p(x^-)=G_x(x,t)+f(t)$, where$x^+=\max (x,0)$,$x^- =\max(-x,0)$,$\phi_p(s)=|s|^{p-2}s$,$p\geq2$, $a $ and $b$ are positive constants $(a\not=b)$ ,the perturbation $f(t)\in {\cal C}^{23}(\RR/2\pi_p \ZZ)$, the oscillating term $G\in {\cal C}^{21}(\RR\times\RR/2\pi_p \ZZ)$,where $\pi_p=\frac{2\pi(p-1)^{\frac{1}{p}}}{p\sin\frac{\pi}{p}},$ and $G(x,t)$ satisfies $\label{G} |D_x^iD_t^jG(x,t)|\le C,\quad 0\le i+j\le 21,$ and $\label{hatG} |D_t^j\hat{G}|\le C,\quad 0\le j\le 21$ for some $C>0$, where $\hat{G}$ is some function satisfying $\frac{\pa \hat{G}}{\pa x}=G$.
\end{abstract}

\vskip 1.0cm
\section{Introduction}








One of the  most studied  semilinear Duffing's equations is
\begin {equation}
\label {tag1.1}
x^{\prime\prime} +ax^+ -bx^- =f(x, t),
\end {equation}
where $x^+=\max (x, 0),\  x^- =\max(-x, 0),\ f(x, t) $ is a smooth
$2\pi$-periodic function on $t$, $a $ and $b$ are positive constants $(a\not=b)$.

If $f(x,t)$ depends only on $t$, the equation (\ref{tag1.1}) becomes
\beq\label{tag1.01}
x^{\prime\prime} +ax^+ -bx^- =f(t),\quad f(t+2\pi)=f(t),
\eeq
which had been studied by Fucik \cite{[Fucik]} and Dancer \cite{[Dancer]} in their
investigations of boundary value problems associated to  equations with
``jumping nonlinearities''. For recent developments, we refer to
\cite{[Gall], Habets, [Lazer]} and references therein.

In 1996, Ortega \cite{[O-London]} proved the Lagrangian stability for the equation
\begin{equation}
\label {tag1.2}
  x'' +ax^+ -bx^- =1+ \gamma h(t)
\end{equation}
if $|\gamma | $ is sufficiently small and $h\in {\cal C}^4(\SS^1)$.

  On the other hand, when $\frac {1}{\sqrt{a}} +\frac {1}{\sqrt{b}}\in \QQ$,
Alonso and Ortega \cite{[Or]} proved that there is a $2\pi$-periodic function
$f(t)$ such that all the solutions of  Eq.~(\ref{tag1.01}) with large
initial conditions are unbounded. Moreover for such a $f(t)$, Eq.~(\ref{tag1.01})
has periodic solutions.

In 1999, Liu \cite{Liujmaa}  removed the smallness assumption on $|\gamma |$ in
  Eq.~(\ref{tag1.2}) when $\frac{1}{\sqrt {a}} +\frac {1}{\sqrt{b}} \in \QQ $ and obtained the same result.

  For the more general equation
\beq
x''+ax^+-bx^-+\phi(x)= e(t)
\eeq
Wang \cite{[Wangx]} and Wang \cite{[Wangy]} considered the Lagrangian stability when the perturbation $\phi(x)$ is bounded.And Yuan \cite{[Yuanx]} investigated the existence of quasiperiodic solutions and Lagrangian stability when  $\phi(x)$ is unbounded.

Fabry and Mawhin  \cite{[FM]} investigated the equation
\beq\label{semi}
x''+ax^+-bx^-=f(x)+g(x)+e(t)
\eeq
under some appropriate conditions, they get the boundedness of all solutions.\\

Yang \cite{[Yangx2]} considered  more complicated nonlinear equation with p-Laplacian operator
\beq\label{semi}
((\phi_p(x'))'+(p-1)[a\phi_p(x^+)-b\phi_p(x^-)]+f(x)+g(x)=e(t).
\eeq

Using Moser's small twist theorem, he proved that all the solutions are bounded,when $\frac{1}{a^{\frac{1}{p}}}+\frac{1}{b^{\frac{1}{p}}}=\frac{2m}{n}$, $m,n\in \NN $, the perturbation $f(x)$ and the oscillating term $g$ are bounded. For the case when $\frac{1}{a^{\frac{1}{p}}}+\frac{1}{b^{\frac{1}{p}}}=2\omega^{-1} $, where $\omega \in \RR^+ \backslash \QQ$, the perturbation $f(x)$ is bounded, Yang \cite{[Yangx1]} studied the following equation
\beq
(\phi_p(x'))'+a\phi_p(x^+)-b\phi_p(x^-)+f(x)=e(t).
\eeq
and came to the conclusion that every solution of the equation is bounded.

In 2004, Liu \cite{[Liujde04]} studied equation
\beq
(\phi_p(x'))'+a\phi_p(x^+)-b\phi_p(x^-)=f(x,t),f(x,t+2\pi)=f(x,t)
\eeq
where $p>1$, for the cases when $\frac{\pi_p}{a^{\frac{1}{p}}}+\frac{\pi_p}{b^{\frac{1}{p}}}=\frac{2\pi}{n}$
and $f\in{\cal C}^{(7,6)}(\RR\times \RR/2\pi\ZZ)$ and satisfies that

   {\rm \ (i)}\  the following limits exists {\rm uniformly\ in\ } $t$
  $$
  \  \lim_{x\rightarrow \infty}f(x,t)=f_{\pm}(t)\quad
$$

{\rm (ii)}\ the following limits exists {\rm uniformly\ in\ } $t$
$$
\ \lim_{x\rightarrow \infty}x^m \frac{\pa^{m+n}}{\pa x^m\pa t^n}f(x,t)=f_{\pm,m,n}(t)
$$
for $(n,m)=(0,6),\ (7,0)$ and $(7, 6)$. Moveover, $f_{\pm,m,n}(t)\equiv 0$ for $m=6,\ n=0, 7$. He comes to the conclusion that all solutions are bounded and the existence of quasi-periodic solutions.\\

In 2012,Jiao,Piao and Wang \cite{[JPW]} considered the bounededness of equations
\beq
x''+\omega^2x+\phi(x)=G_x(x,t)+f(t),
\eeq and
\beq
x''+ax^+-bx^-=G_x(x,t)+f(t).
\eeq

Inspired by the above references, we are going to study the boundedness of all solutions for the more general equation
\beq
\label{tag(1.10)}
(\phi_p(x'))'+a\phi_p(x^+)-b\phi_p(x^-)=G_x(x,t)+f(t)
\eeq
 Our main results are as follows:\\

\begin{Theorem}
   Assume $f(t)\in {\cal C}^{23}(\RR/2\pi_p \ZZ)$, $G\in {\cal C}^{21}(\RR\times\RR/2\pi_p \ZZ)$,and $G(x,t)$ satisfies
\beq\label{G} |D_x^iD_t^jG(x,t)|\le C,\quad 0\le i+j\le 21\eeq and \beq\label{hatG} |D_t^j\hat{G}|\le C,\quad 0\le j\le 21\eeq for some $C>0$, where $\hat{G}$ is some function
satisfying $\frac{\pa \hat{G}}{\pa x}=G$,and $\omega=\frac{1}{2}(\frac{1}{a^{\frac{1}{p}}}+\frac{1}{b^{\frac{1}{p}}})$, and $\omega \in \RR^+ \backslash \QQ$
satisfy the Diophantine condition:
\begin{equation}\label{Diophantine}|m\omega +n|\ge
\frac{\gamma}{|m|^{\tau}},\quad \forall \ \ (m, n)\not=(0,0)\in
{\Bbb Z}^{2}, \end{equation} where $1<\tau<2,\ \gamma>0$, and $[f]=\frac{1}{2\pi_p}\int_0^{2\pi_p} f(t) dt \not= 0, $ where $\pi_p=\frac{2\pi(p-1)^{\frac{1}{p}}}{p\sin\frac{\pi}{p}}$.
Then equation (\ref{tag(1.10)}) possesses Lagrange stability, i.e. if $x(t)$ is any solution of equation (\ref{tag(1.10)}), then it exists for all $t\in R$ and $\sup_{t\in R}(|x(t)|+|\dot{x}(t)|)<\infty$.
\end{Theorem}
\begin{Remark}
In the above, $\gamma$ can be any positive number. Thus our statement holds true for $\omega$ of full measure.
\end{Remark}
\begin{Remark}
In Liu[17], it is required that  f satisfies the limit condition, which is not satisfied by the function
G in our situation. Thus our situation is more general.
\end{Remark}

The main idea is as follows: By means of transformation theory the original system outside of a
large disc $D=\{(x, x')\in R^{2}:  x^{2}+x'^{2}\leq r^{2} \}$ in
$(x, x')$-plane is transformed into a perturbation of an integrable
Hamiltonian system. The Poincar\'{e} map of the transformed system
is closed to a so-called twist map in $R^{2}\backslash D$. Then
Moser's twist theorem guarantees the existence of arbitrarily large
invariant curves diffeomorphic to circles and surrounding the origin
in the $(x,x')$-plane. Every such curve is the base of a
time-periodic and flow-invariant cylinder in the extended phase
space $(x,x',t)\in R^{2}\times R,$ which confines the solutions in
the interior and which leads to a bound of these solutions.\\

The remain part of this paper is organized as follows. In section 2,
we introduce action-angle variables and exchange the role of time and angle variables. In section 3, we construct
canonical transformations such that the new Hamiltonian system is
closed to an integrable one. In section 4, we will prove the
Theorem 1 by Moser's twist theorem.

Throughout this paper, $F(x)=\int_0^x f(s)ds,F(0)=0$, $c$ and $C$  are some positive constants without concerning their quantity.


\section{Some Canonical transformations}
In this section,  we will state  some technical lemmas which will be
used in the proof of Theorem 1. Throughout this section,  we
assume the hypotheses of Theorem 1 hold.

\vskip 2mm
\subsection{Action-angle variables}
Borrowing the idea from Liu \cite {[Liujde04]} and Yang \cite {[Yangx1]},we introduce a new variables $y$ as $y=-\varphi_p(\omega x)$, let $q$ be the conjugate exponent of $p$ : $p^{-1}+q^{-1}=1$. Then (\ref{tag(1.10)}) is changed into the form
\beq
x'=-\omega^{-1}\varphi_q(y),   y'=\omega^{-1}[a_1\varphi_p(x^{+})-b_1\varphi_p(x^{-})]-\omega^{p-1}[G_x(x,t)+f(t)]
\eeq
where $a=\omega^{-p}a_1$, $b=\omega^{-p}b_1$ and $a_1$, $b_1$ satisfy
\beq
a_1^{-\frac{1}{p}}+b_1^{-\frac{1}{p}}=2,
\eeq
which is a planar non-autonomous Hamiltonian system
\beq
x'=-\frac{\pa H}{\pa y }(x,y,t),  y'=\frac{\pa H}{\pa y }(x,y,t)
\eeq
where $$H(x,y,t)=\frac{\omega^{-1}}{q}|y|^q+\frac{\omega^{-1}}{p}(a_1|x^+|^p+b_1|x^-|^p)-\omega^{p-1}[G(x,t)+f(t)x]. $$
Let $C(t)=\sin_p t$ be the solution of the following initial value problem
\beq
(\varphi_p(C'(t)))'+\varphi_p(C(t))=0,
\quad C(0)=0,C'(0)=1.
\eeq
  Then it follows from [16] that $C(t)=\sin_p(t)$ is a $2\pi_p$-period $C^2$ odd function with $\sin_p(\pi_p-t)=\sin_p(t)$, for $t\in[0,\frac{\pi_p}{2}]$ and
$\sin_p(2\pi_p-t)=-\sin_p(t)$, for $t\in[\pi_p,2\pi_p]$. Moreover for $t\in(0,\frac{\pi_p}{2})$, $C(t)>0$, $C'(t)>0$, and $C: [0,\frac{\pi_p}{2}]\rightarrow[0,(p-1)^{\frac{1}{p}}]$ can be implicitly given by
$$\int_0^{\sin_p t}\frac{ds}{(1-\frac{s^p}{p-1})^{\frac{1}{p}}}=t. $$
\begin{Lemma}
For $p\geq2$ and for any $(x_0,y_0)\in R^2$, $t_0\in R$, the solution $$z(t)=(x(t, t_0, x_0, y_0),y(t, t_0, x_0, y_0))$$ of (2.1) satisfying the initial condition $z(t_0)=(x_0,y_0)$ is unique and exists on the whole $t$-axis.
\end{Lemma}
 The proof of uniqueness can be obtained similarly as the proof of Proposition 2 in [17], the global existence result can be proved similarly as Lemma 3.1 in [10].
Consider an auxiliary equation
\[(\phi_p(x'))'+a_1\phi_p(x^+)-b_1\phi_p(x^-)=0\]
Let $v(t)$ be the solution with initial condition: $(v(0),v'(0))=((p-1)^{\frac{1}{p}},0)$. Setting $\phi_p(v')=u $, then $(v,u)$ is a solution of the following planar system:
\[x'=\phi_q(y),\quad y'=-a_1\phi_p(x^+)+b_1\phi_p(x^-)\]
where $q=p/(p-1)>1$. It is not difficult to prove that:\\
(i) $q^{-1}|u|^q+p^{-1}(a_1|v^+|^p+b_1|v^-|^p)\equiv\frac{a_1}{q}$;\\
(ii) $v(t)$ and $u(t)$ are $2\pi_p$-periodic functions.\\
(iii)$v(t)$ can be given by
\beq
v(t)=\left\{
\begin{array}{ll}
\sin_p(a_1^{\frac{1}{p}}t+\frac{\pi_p}{2}),&\mbox{$0\leq t\leq\frac{\pi_p}{2a_1^{\frac{1}{p}}}$,}\\
-(\frac{a_1}{b_1})^{\frac{1}{p}}\sin_p{b_1}^{\frac{1}{p}}(t-\frac{\pi_p}{2a_1^{\frac{1}{p}}}),& \mbox{$
\frac{\pi_p}{2a_1^{\frac{1}{p}}}<t\leq\pi_p$.}
\end{array}
\right.
\eeq
\beq
v(2\pi_p-t)=v(t), t\in[\pi_p,2\pi_p].
\eeq
\begin{Lemma}
Let $I_p=\int_0^{\frac{\pi_p}{2}}\sin_p tdt.$ Then $$I_p=\frac{(p-1)^{\frac{2}{p}}}{p}B(\frac{2}{p},1-\frac{1}{p}),$$
where $B(r,s)=\int_0^1t^{r-1}(1-t)^{s-1}dt$ for $r>0,s>0.$
\end{Lemma}
 From the expression of $v(t)$ in (2.5), we obtain
\beq
\int_0^{\frac{\pi_p}{2a_1^{\frac{1}{p}}}}v(t)dt=\frac{I_p}{a_1^{\frac{1}{p}}},
\eeq
\beq
\int_{\frac{\pi_p}{2a_1^{\frac{1}{p}}}}^{\pi_p}v(t)dt=-\frac{a_1^{\frac{1}{p}}I_p}{b_1^{\frac{2}{p}}}.
\eeq
This method has been used in [8].

We introduce the action and angle variables via the solution $(v(t),u(t))$ as follows.
\[x=d^{\frac{1}{p}}r^{\frac{1}{p}}v(\theta),y=d^{\frac{1}{q}}r^{\frac{1}{q}}u(\theta)\]
where $d=pa_1^{-1}$. This transformation is called a generalized symplectic transformation as its Jacobian is 1. Under this transformation, the system (2.1) is changed to
\beq\label{h}
\theta'=\frac{\partial h}{\partial r}(r,\theta,t),r'=-\frac{\partial h}{\partial \theta}(r,\theta,t)
\eeq
with the Hamiltonian function
\beq\label{Hamiltonian2}
h(r,\theta,t)=\omega^{-1} r-f_1(r,\theta,t)-\omega^{p-1}d^{\frac{1}{p}}r^{\frac{1}{p}}v(\theta)f(t)
\eeq
where $f_1(r,\theta,t)=\omega^{p-1}G( d^{\frac{1}{p}}r^{\frac{1}{p}}v(\theta),t).$\\

For any function $f(\cdot, \t)$, we denote by $[f](\cdot)$ the
average value of $f(\cdot, \theta) $ over $\mathbb{S}_p\triangleq\RR/2\pi_p \ZZ$, that is,
$$
[f](\cdot):=\frac {1}{2\pi_p} \int _0^{2\pi_p} f (\cdot, \theta) d\theta.
$$\\

For the above function $f_1(r,\theta,t)$ in (2.10) we have
\begin{Lemma}\label{lemmaf1}
The following conclusion holds true: \beq\label{f1}
|D_r^iD_t^jf_1(r,\theta,t)|\le C\cdot r^{-\frac{i}{q}},\quad 0\le
i+j\le 21. \eeq
\end{Lemma}
\Proof
The proof of this lemma can get directly from the definition of $f_1$ and the conditions in Theorem 1.

\noindent The following technique lemma will be used to refine the estimates on $[f_1](r, t)$.
\begin{Lemma}
 Assume $f\in {\cal C}^1(\RR/2\pi_p\ZZ)$, $G(x,t)\in {\cal C}^1(\RR^1\times \RR/2\pi_p\ZZ)$ and
 $G_x'(x,t)=g(x,t)$. Suppose there are two positive constants $\bar{G}$ and $\bar{g}$ such that $|G(x,t)|\le\bar{G}$, $|g(x,t)|\le \bar{g}$ for any $(x,t)$.  Let $A(r,\t)\in {\cal C}^2(\RR^1\times \RR/2\pi_p\ZZ)$ be of the form $A(r, \t)=(r+h(r,\t))^{\frac{1}{p}}$ with
 \beq \label{h1}
 h, \frac{\pa h}{\pa \t}, \frac{\pa^2 h}{\pa \t^2}=O(r^{\frac{1}{p}})
  \eeq
  for $r\gg 1.$

  Then for any
 constant $\delta_0\in (0, \frac{1}{10})$ it holds that
 \beq\label{partialintegral}
\left |\int_0^{2\pi_p}f(\theta)g(A v(\theta),t)d\t\right |\le C\cdot
r^{-\delta_0},\quad r\gg 1,
 \eeq
 where $C$ depends only on $\bar{G}$, $\bar{g}$ and $\|f\|_{C^0}$.
\end{Lemma}
\Proof
 Let $[0, 2\pi_p]=I_1\bigcup I_2$,  where
 $I_1=[0,r^{-2\delta_0}]\bigcup[\pi_p-r^{-2\delta_0},
 \pi_p+r^{-2\delta_0}]\bigcup[2\pi_p-r^{-2\delta_0}, 2\pi_p]$ and
 $I_2=[r^{-2\delta_0}, \pi_p-r^{-2\delta_0}]\bigcup[\pi_p+r^{-2\delta_0},
 2\pi_p-r^{-2\delta_0}]$. Then
 $$
\int_0^{2\pi}f(\theta)g(A v(\theta),t)d\t=\int_{I_1}f(\theta)g(A v(\theta),t)d\t+\int_{I_2}f(\theta)g(A v(\theta),t)d\t.$$ Obviously, $|I_1|\le
C\cdot r^{-2\delta_0}$, where $|\cdot|$ denotes the Lesbegue measure.
Then from the boundedness of $g(x,t)$, it is easy to see that
$$ \left|\int_{I_1}f(\theta)g(A v(\theta),t)d\t\right|\le C\cdot r^{-2\delta_0}.$$
To estimate the integral on $I_2$, we first estimate the integral on
the interval $I_{21}=[r^{-2\delta_0}, \pi_p-r^{-2\delta_0}]$.

Consider
$D_{\t}(A v(\theta))=A'_{\t}v(\theta)-A v'(\theta).$ From (\ref{h1}), it holds that $|Av'(\theta)|\ge c\cdot r^{\frac{1}{p}-2\delta_0}$ and $A'_{\t}\cdot v(\theta)=O(1)$ for $\t\in I_{21}$, which implies
\beq \label{DA}|D_{\t}(A v(\theta))|\ge c\cdot r^{\frac{1}{p}-2\delta_0}. \eeq

Similarly from the definition of $A$ and the condition (\ref{h1}), we have
\beq\label{D2A}D^2_{\t}(A v(\theta))=D^2_{\t}A\cdot v(\theta)+2 D_{\t}A\cdot v'(\theta)+A v''(\theta)=O(r^{\frac{1}{p}}). \eeq
By direct computation, we have
$$D_{\t}(f(\theta)(D_{\t}(A v(\theta)))^{-1})=f'\cdot (D_{\t}(Av(\theta)))^{-1}+f\cdot (D_{\t}(Av(\theta)))^{-2}\cdot (-D^2_{\t}(Av(\theta))).$$ Thus from (\ref{DA}) and (\ref{D2A}), we obtain the estimate
\beq\label{antiderivative}|D_{\t}(f(\theta)(D_{\t}(A v(\theta)))^{-1})|\le C\cdot r^{4\delta-\frac{1}{p}}.\eeq

By
integration by parts, we have that
$$\begin{array}{ll}
&\int_{I_{21}}f(\theta)g(Av(\theta),t)d\t=\int_{I_{21}}f(\theta)(D_{\t}(Av(\theta)))^{-1}dG(Av(\theta),t)\\
\\
=&(D_{\t}(Av(\theta)))^{-1}
f(\t)G(Av(\theta),t)|_{r^{-2\delta_0}}^{\pi_p-r^{-2\delta_0}}-\int_{I_{21}}
G(Av(\theta),t) D_{\t}(f(\theta)D_{\t}((Av(\theta)))^{-1})d\t.
\end{array}
$$
From (\ref{DA}) and (\ref{antiderivative}) , for $\t\in I_{21}$ it holds that
$$
\left|(D_{\t}(Av(\theta)))^{-1} f(\t)G(Av(\theta),t)|_{\t=r^{-2\delta_0}}\right
|,\ \left |(D_{\t}(Av(\theta)))^{-1}
f(\t)G(Av(\theta),t)|_{\t=\pi_p -r^{-2\delta_0}}\right|\le C\cdot
r^{4\delta_0-\frac{1}{p}}
$$
and
$$
\left|G(Av(\theta),t)\cdot D_{\t}(f(\theta)D_{\t}((Av(\theta)))^{-1})\right|\le C\cdot r^{4\delta_0-\frac{1}{p}}.\quad
$$

 Similarly, we can have the same estimate for the
other parts of $I_2$.

Hence from the fact $0<\delta_0<\frac{1}{10}$, we obtain
(2.7). The proof of this lemma is completed.
\qedbox
 \vskip 0.4cm
For $[f_1](r,t)$, we have the following
result:
\begin{Corollary}
\label{refinelemmaf1} The following conclusion holds true:
\beq\label{refinef1} |D_r^iD_t^j[f_1](r,t)|\le C\cdot
r^{-\delta_1-\frac{i}{p}},\quad 0\le i+j\le 21,\eeq where the
constant $\delta_1$ is in $(0,\frac{1}{10})$.
\end{Corollary}
\Proof From the definition of $f_1$, we have
$[f_1](r,t)=\frac{1}{2\pi_p}\int_0^{2\pi_p}G(r^{\frac{1}{p}}v(\theta),t)d\t$.
From (\ref{G}) and (\ref{hatG}), we know that $G$ and $\hat{G}$ are
bounded. Thus for $i+j=0$, (\ref{refinef1}) is deduced from lemma
2.2 where we set $f\equiv 1$ and $A(r,\t)=r^{\frac{1}{p}}$. For
$i+j\ge 1$, it can be easily seen that $\frac{\pa^{i+j}}{\pa
r^{i}\pa t^j}G$ are the sum of the term like
$$
\frac{\pa^{k+j}}{\pa x^{k}\pa
t^j}G(r^{\frac1p}v(\theta),t)(r^{\frac{1}{p}})^{(i_1)}\cdots
(r^{\frac{1}{p}})^{(i_k)}\cdot (v(\theta))^k,
$$
where $i_1+\cdots i_k=i$. Thus (\ref{refinef1}) is implied from
lemma 2.2 for the function $\frac{\pa^{k+j}}{\pa x^{k}\pa
t^j}G(r^{\frac{1}{p}}v(\theta),t)$ and (\ref{G}). This ends the proof of
the lemma. \qedbox

\vskip 3mm \subsection{Exchange of the roles of time and angle
variables}

According to Levi \cite{LE1},the equality
$$
rd\theta -h dt =-(hdt -rd\theta ),
$$
means if we can solve $r=r(h,t,\theta) $ from
Eq.(\ref{h}) as a function of $h, t $ and $\theta $, then we have
\begin {equation}
\label {5.1} \frac{dh}{d\theta}=-\frac{\partial r}{\partial t}(h, t,
\theta), \ \ \frac {dt}{d\theta}=\frac{\partial r}{\partial
h}(h,t,\theta),
\end {equation}
i.e., Eq.(\ref{5.1}) is a Hamiltonian system with Hamiltonian
function $r=r(h,t, \theta) $ and now the action, angle and time
variables are $h, t, $ and $\theta $,  respectively.

From Eq.(\ref{Hamiltonian2}) and lemmas, it follows that
$$
\lim_{r\rightarrow +\infty}\frac{h}{r} =\omega^{-1} > 0
$$
and for $r\gg 1$
$$
\frac { \partial h}{\partial r}=\omega^{-1}
-\frac{\partial}{\partial r} f_1 (r,
\theta)-\frac{1}{p}f(t)\omega^{p-1}d^{\frac{1}{p}}r^{\frac{1}{p}-1}v(\theta)>0.
$$
By the implicit function theorem,  we know that there is a function
$R=R(h, t, \theta) $ such that
\begin {equation}
\label{r} r(h, t, \theta)=\omega h -R (h,t,\theta).
\end {equation}
Moreover,  for $h\gg 1$,
$$
|R(h,t,\theta)| \le \omega h /2
$$
and $R(h, t, \theta) $ is $C^{19} $ in $h $ and $ t$.



From (\ref{Hamiltonian2}), it holds that \beq\label{RR}
R=\omega f_1(\omega h-R,
t,\theta)-\omega^p d^{\frac{1}{p}}{(\omega h-R)}^{\frac{1}{p}}v(\theta)f(t). \eeq
The proof of following two lemmas are slightly different to \cite{Liu1999}, here for the convenience of readers, we give the proofs of them.
\begin{Lemma}\label{lemmaR}
Assume $R$ is defined by (\ref{RR}) with $|R|\ll h$ for $h\gg 1$.
 Then
it holds that \beq\label{prop-R} |D_h^iD_t^jR|\le C\cdot
h^{n(i)},\quad 0\le i+j\le 21 \eeq for $h\gg 1$, where $n(i)=-\frac{i}{q}$
for $i\ge 1$ and $n(0)=\frac{1}{p}$.
\end{Lemma}
Proof. \noindent{\rm (i)}\ $i+j=0.$\quad The proof for this case can be easily obtained from lemma \ref{lemmaf1} and the conditions in the Theorem .

\noindent{\rm (ii)}\ $i+j=1.$\quad It is clear that for $h\gg 1$,
$$
|\o\frac{\pa f_1}{\pa r}(\o h-R,t,\t)|+|{\frac{\o^p}{p}}d^{\frac{1}{p}}(\o h-R)^{-\frac{1}{q}}v(\theta) f(t)|\le \frac{1}{2}.
$$
Define
$$
\Delta(h,t,\t)=1+\o\frac{\pa f_1}{\pa
r}(\o h-R,t,\t)-\frac{\o^p}{p}d^{\frac{1}{p}}(\o h-R)^{-\frac{1}{q}}v(\theta) f(t),$$
$$
g_1=\o^{2}\frac{\pa f_1}{\pa r}(\o h-R,t,\t)-\frac{\o^{p+1}}{p}d^{\frac{1}{p}}(\o h-R)^{-\frac{1}{q}}v(\theta) f(t),$$
$$
g_2=-\omega^p d^{\frac{1}{p}}(\o h-R)^{\frac{1}{p}}v(\theta) f(t)+\o\frac{\pa f_1}{\pa t}(\o h-R,t,\t).$$
Then it follows that
\beq\label{g_1g_2}
\Delta \cdot \frac{\pa R}{\pa h}=g_1,\quad \Delta\cdot \frac{\pa R}{\pa t}=g_2.
\eeq

From lemma \ref{lemmaf1},$p\geq2$ and the boundedness of $f(t)$, we have $|g_1|\le C\cdot h^{-\frac{1}{q}}$
and $|g_2|\le C\cdot h^{\frac{1}{p}}$. Thus the proof for this case
is completed.

{\rm(iii)} $i+j=2$.\quad Lemma \ref{lemmaf1} implies that
$$
|\frac{\pa \Delta}{\pa t}|\le C\cdot h^{-\frac{1}{q}},\ |\frac{\pa \Delta}{\pa h}|\le C\cdot h^{-\frac{2}{q}},\ |\frac{\pa g_1}{\pa t}|\le C\cdot h^{-\frac{1}{q}},\ |\frac{\pa g_1}{\pa h}|\le C\cdot h^{-\frac{2}{q}},
\ |\frac{\pa g_2}{\pa h}|\le C\cdot h^{-\frac{1}{q}},\ |\frac{\pa g_2}{\pa t}|\le C\cdot h^{\frac{1}{p}}.$$
From the second equation of (\ref{g_1g_2}), we obtain
$$
\Delta \frac{\pa^2 R}{\pa t^2}+\frac{\pa \Delta}{\pa t}\cdot \frac{\pa R}{\pa t}=\frac{\pa g_2}{\pa t}$$
and
$$
\Delta \frac{\pa^2 R}{\pa t\pa h}+\frac{\pa \Delta}{\pa h}\cdot \frac{\pa R}{\pa t}=\frac{\pa g_2}{\pa h}.$$
The above inequalities and equations imply that
$$
|\frac{\pa^2 R}{\pa t^2}|\le C\cdot h^{\frac{1}{p}},\quad |\frac{\pa^2R}{\pa h\pa t}|\le C\cdot h^{-\frac{1}{q}}.
$$
From the first equation of (\ref{g_1g_2}), we know that
$$
\Delta\frac{\pa^2R}{\pa h^2}+\frac{\pa \Delta}{\pa h}\cdot \frac{\pa R}{\pa h}=\frac{\pa g_1}{\pa h},
$$
which implies $|\frac{\pa^2 R}{\pa h^2}|\le C\cdot h^{-\frac{2}{q}}$. Thus we complete the proof for this case.
\vskip 0.3cm
In general, if
$$
|D_h^iD_t^jR|\le C\cdot h^{n(i)},\quad 0\le i+j\le m,
$$
then it holds that
$$
|D_h^iD_t^j\Delta|\le C\cdot h^{-\frac{1}{q}+n(i)},\quad |D_h^iD_t^jg_1|\le C\cdot h^{-\frac{1}{q}-\frac{i}{q}},\quad
|D_h^iD_t^jg_2|\le C\cdot h^{-\frac{i}{q}}.
$$
Consequently, we obtain
$$
|D_h^iD_t^jR|\le C\cdot h^{n(i)},\quad 0\le i+j\le m+1.
$$
The proof is completed.
\qedbox
\vskip 0.3cm

In (\ref{RR}), we denote $R=-\omega^p d^{\frac{1}{p}}{(\omega h)}^{\frac{1}{p}}v(\theta)f(t)-R_1(h,t,\t)$. Then
\beq\label{R_1}
R_1=\omega f_1(\omega h-R,t,\theta)-{\frac{1}{p}}\int_0^1\omega^p d^{\frac{1}{p}}{(\o h-\tau R)}^{-\frac{1}{q}}Rv(\theta)f(t)d\tau.
\eeq
Then we have the following conclusion:
\begin{Lemma}\label{lemmaR1}
It holds that
$$|D_h^iD_t^jR_1|\le C\cdot h^{-\frac{i}{q}},\quad 0\le i+j\le 21.
$$
\end{Lemma}
Proof. The lemma is easily followed from the following claim:
\vskip 0.3cm
{\bf Claim}
\beq\label{claim}
\begin{array}{ll}
&|D_h^iD_t^jf_1(\o h-\tau R,t,\t)|\le C\cdot h^{-\frac{i}{q}},\\
&|D_h^iD_t^j(\o h-\tau R)^{-\frac{1}{q}}d^{\frac{1}{p}}Rv(\theta)f(t)|\le C\cdot h^{-\frac{1}{q}-\frac{i}{q}}
\end{array}
\eeq
for $0\le i+j\le 21$.
\vskip 0.3cm
{\it Proof of the claim.}\quad
We only prove the first inequality of (\ref{claim}) and the proof for the other is similar.

{\rm (i)}\ $i+j=0$. The proof for this case can be obtained directly
from lemma 2.1.

{\rm (ii)}\ $i>0, j=0$. We have the following equality:
$$
D_h^if_1(\o h-\tau R,t,\t)=\sum \frac{\pa^k f_1}{\pa r^k}(u,t,\t)\cdot \frac{\pa^{i_1} u}{\pa h^{i_1}}\cdots
\frac{\pa^{i_k} u}{\pa h^{i_k}}
$$
with $0< k\le i,\ i_1,\cdots, i_k>0,\ i_1+\cdots i_k=i$ and $u=\o h-\tau R$. Assume there are $l$($\le k$) numbers in $\{i_1,\cdots, i_k\}$
which is equal to 1. Then we obtain
$$
|D_h^i f_1(u,t,\t)|\le C\cdot h^{-\frac{k}{q}}\cdot h^{-\frac{i_1+\cdots i_k-l}{q}}\le C\cdot  h^{-\frac{i}{q}}.
$$
{\rm (iii)}\ $i=0, j>0$. By direct computation, we have
$$
D_t^jf_1(\o h-\tau R,t,\t)=\sum \frac{\pa^{k+l}f_1}{\pa r^{k}\pa t^l}(u,t,\t)\cdot \frac{\pa^{j_1} u}{\pa t^{j_1}}\cdots
\frac{\pa^{j_k} u}{\pa t^{j_k}}
$$
with $0\leq k\le j,0\leq l\le j,k+l=j,\ j_1,\cdots, j_k>0,\ j_1+\cdots j_k=k$. It follows that
$$
|D_t^jf_1(u,t,\t)|\le C\cdot h^{-\frac{k}{q}}\cdot h^{\frac{k}{p}}\le C.
$$
The last step,we get from that $p\geq2$, $\frac{1}{p}+\frac{1}{q}=1$ ,and $\frac{1}{p}\leq\frac{1}{q}$.

{\rm (iv)}\ $i>0, j>0$. By direct computation, we have
$$
D_h^iD_t^j \frac{\pa f_1}{\pa r}(u,\t)=\sum \frac{\pa^{k_1+k_2+l}f_1}{\pa r^{k_1+k_2} \pa t^l}(u,\t)\cdot\frac{\pa^{i_1} u}{\pa h^{i_1}}\cdots
\frac{\pa^{i_{k_1}} u}{\pa h^{i_{k_1}}}\cdot \frac{\pa^{l_1+j_1} u}{\pa h^{l_1}\pa t^{j_1}}\cdots
\frac{\pa^{l_{k_2}+j_{k_2}} u}{\pa h^{l_{k_2}}\pa t^{j_{k_2}}},
$$
where $u=\o h-\tau R$ and
$$
0\le k_1\le i,\ 0\le k_2\le j,\ 0\leq l\leq j,\ k_2+l=j,\ i_1,\cdots,i_{k_1},\ j_1,\cdots, j_{k_2}>0,\quad l_1,\cdots, l_{k_2}\ge 0,
$$
$$
i_1+\cdots i_{k_1}+l_1+\cdots+ l_{k_2}=i,\quad j_1+\cdots+j_{k_2}+l=j.
$$
Assume that there are $m$($\le k_1$) numbers in $\{i_1,\cdots, i_{k_1}\}$ which is equal to 1. Then
$$
|D_h^iD_t^j \frac{\pa f_1}{\pa r}|\le C\cdot h^{-\frac{k_1+k_2}{q}}\cdot h^{-\frac{i_1+\cdots+i_{k_1}+l_1+\cdots+l_{k_2}-m}{q}}\le C\cdot
h^{-\frac{i}{q}}.
$$
This ends the proof of the claim.
\qedbox
\vskip 0.3cm
From the definition of $R_1$, we can obtain the following conclusion:
\begin{Lemma}
\label{[R_1]}
For the function $[R_1](h, t)$, we have that
$$|D_h^iD_t^j[R_1]|\le C\cdot (h^{-i}+h^{-\delta_1-\frac{i}{q}}),\quad 0\le i+j\le 21,
$$
where $\delta_1\in (0, \frac{1}{10})$.
\end{Lemma}

\vskip 0.3cm

From (\ref{r}), (\ref{RR}) , we obtain that the Hamiltonian $r(h,t,\t)$ in (\ref{r}) is of the form:
 \beq\label{Hamiltonianr} r=\omega h+\omega^p d^{\frac{1}{p}}{(\omega h)}^{\frac{1}{p}}v(\theta)f(t)+R_1(h,t,\t). \eeq


\section{More canonical transformations}
In this section, we will make some more canonical transformations
such that the Poincar\'e map of the new system is close to twist
map.

\begin{Lemma}\label{lemma3.1}
There exists a canonical transformation $\Phi_1$ of the form:
$$\Phi_1:\quad \left\{\begin{array}{ll}
h&=\rho\\
t&=\tau+V_1(\rho,\tau,\t)
\end{array}\right.
$$ where the functions $ V_1$ are periodic in $\tau,
\theta$. Under this transformation, the Hamiltonian system with
Hamiltonian (\ref{Hamiltonianr}) is changed into the following
one \beq \label{tilder} \tilde{r}=\omega
\rho+\omega^p d^{\frac{1}{p}}(\omega\rho)^{\frac{1}{p}}v( \theta)[f]+\tilde{R}_1(\rho,\tau,\theta),
\eeq
Moreover, the new perturbation $\tilde{R}_1$ satisfies
\beq\label{[tilder1]}|\frac{\partial ^{i+j}}{\partial \rho ^i \partial \tau^j}\tilde{R}_1|\le C\cdot
\rho^{-\frac iq},\quad 0\le i+j\le 21.\eeq Moreover, for the function$[\tilde{R}_1](\rho,\theta)$,it holds that
\beq\label{[tilder2]}
|D_\rho^iD_\tau^j[\tilde{R}_1]|\le C\cdot (\rho^{-i}+\rho^{-\delta_1-\frac{i}{q}}).
\quad 0\le i+j\le 21.
\eeq
\end{Lemma}
\Proof We construct the canonical transformation by means of
generating function:
$$
\Phi_1:\quad h=\rho,\quad t=\tau+\frac{\pa S_1}{\pa \rho}(\rho,\tau,\t).
$$
Under this transformation,  the new Hamiltonian function $ \tilde{r} $
is of the form
$$
\tilde{r}=\omega\rho+\omega^p d^{\frac{1}{p}}{(\omega\rho)}^{\frac{1}{p}}v(\theta)f(\tau+{\frac{\pa S_1}{\pa \rho}} )+R_1(\rho,\tau+{\frac{\pa S_1}{\pa \rho}},\theta)+{\frac{\pa S_1}{\pa \theta}}
$$
Let $S_1=-\int_0^{\theta}\omega^p d^{\frac{1}{p}}{(\omega\rho)}^{\frac{1}{p}}v(\vartheta){f(t)-[f]}d\vartheta$,then we have
$$
\tilde{r}(\rho,\tau,\theta)=\omega\rho+\omega^p d^{\frac{1}{p}}{(\omega\rho)}^{\frac{1}{p}}v(\theta)[f]+\tilde{R}_1(\rho,\tau,\theta)
$$
where $\tilde{R}_1(\rho,\tau,\theta)=R_1(\rho,\tau+{\frac{\pa S_1}{\pa \rho}},\theta)=R_1(\rho,\tau,\theta)+\int_0^1{\frac{\pa R_1}{\pa t}}(\rho,\tau+s{\frac{\pa S_1}{\pa \rho}},\theta){\frac{\pa S_1}{\pa \rho}}ds$
From {\ref{lemmaR1}} and the definition of $\tilde{R}_1$, we can get the estimates ({\ref{[tilder1]}}) ,({\ref{[tilder2]}}) can get from {\ref{[R_1]}} and the definition of $\tilde{R}_1.$ \qedbox
 \vskip 1cm
\begin{Lemma}\label{lemma3.2}
There exists a
canonical transformation $\Phi_2$ of the form:
$$\Phi_2:\quad \left\{\begin{array}{ll}
\rho&=I\\
\tau&=s+V_2(I,\t)
\end{array}\right.
$$
with $\tilde{T}(I,\theta+2\pi_p)=\tilde{T}(I,\theta)$,such that
 the
system with Hamiltonian (\ref{tilder}) is transformed into the form:
\begin {equation}
\label {3.13}
\frac{\pa I}{\pa \theta}=-\frac{\pa {\bar r}}{\pa s}(I,s,\theta),\ \ \
\frac{\pa s}{\pa \theta}=\frac{\pa {\bar r}}{\pa I}(I,s,\theta)
\end {equation}
with $\bar r(I,s,\theta)=\omega I+c^*I^
{\frac{1}{p}}+\tilde{R}_2(I,s,\theta)$ and
$c^*\not=0$, where we use the fact that $[f]\not=0$.
Moreover, the new perturbation $\tilde{R}_2$ satisfies
 \beq\label{tilder R_2}
|D_{I}^iD_{s}^j\tilde{R}_2)|\le C\cdot
I^{-\frac{i}{q}}, \quad 0\le i+j\le 21. \eeq Moreover, for the function
$[\tilde{R}_2]_0(I)=(\frac{1}{2\pi_p})^2\int_0^{2\pi_p}\int_0^{2\pi_p}\tilde{R}_2(I,s,\theta)dsd\t $, it holds that
 \beq\label{[tilder R_2]} |D_{I}^i[\tilde{R}_2]_0|\le
C\cdot (I^{-i}+I^{-\delta_1-\frac{i}{q}}),\quad 0\le i\le 21. \eeq

\end{Lemma}
\Proof
 The proof is similar to \cite{Liuwang}, but for the convenience of readers we still give a
 detailed argument.
  We shall look for the required transformation $\Phi_2 $ by
means of a generating function $S_2(I, s,\theta)$, so that
$\Phi _2 $ is implicitly defined by \beq\label{Phi} \Phi _2:\quad
\rho=I+\frac{\partial }{\partial s}S_2(I, s,\theta), \quad
\tau=s +\frac{\partial }{\partial I}S_2(I, s,\theta). \eeq
Under this transformation, the system is changed into the form:
$$
\frac{\pa I}{\pa \theta}=-\frac{\pa {\bar r}}{\pa s}(I,s,\theta),\ \ \
\frac{\pa s}{\pa \theta}=\frac{\pa {\bar r}}{\pa I}(I,s,\theta)
$$
the new Hamiltonian function $ \bar{r} $
is of the form
$$
\bar{r}=\omega\rho+\omega^{p+{\frac{1}{p}}}d^{\frac{1}{p}}[f]\rho^{\frac{1}{p}}v(\theta)+\tilde{R_1}(\rho,\tau,\theta)+{\frac{\pa S_2}{\pa \theta}}
$$
Now we choose
$$
S_2=-\int_0^{\theta}\omega^{p+{\frac{1}{p}}}d^{\frac{1}{p}}[f]\rho^{\frac{1}{p}}v(\vartheta)-c^{*}\rho^{\frac{1}{p}}d\vartheta
$$
where
$
c^{*}=\omega^{p+{\frac{1}{p}}}d^{\frac{1}{p}}[f]\neq0.
$ Obviously,$S_2$ does not depend on $s$ and it is $2\pi_p-$periodic in ${\theta}.$ Hence $\rho=I.$
Let $$\tilde{T}(I,\theta)={\frac{\pa S_2}{\pa I}}.
$$
Then the canonical transformation $\Phi_2 $ is of the form
$$
\rho=I, \tau=s+\tilde{T}(I,\theta).
$$\
Let
\beq
\tilde{R}_2(I,s,\theta)=\tilde{R}_1(\rho,s,\theta)+\int_0^1{\frac{\pa \tilde{R}_1}{\pa \tau}}(\rho,s+m\tilde{T},\theta)\tilde{T}d m.
\eeq
From ({\ref{[tilder1]}}) in lemma{\ref{lemma3.1}},we can get ({\ref{tilder R_2}}) easily, and ({\ref{[tilder R_2]}) can get from ({\ref{[tilder2]}}).
The proof of this lemma is completed.
\quad \qedbox\

For convenience, we denote
\beq\label{barr}
\bar{r}=\omega I+\bar{r}_1(I)+\bar{r}_2(I,s,\theta),
\eeq
with $\bar{r}_1=c^{*}I^{\frac{1}{p}},$ $\bar{r}_2(I,s,\theta)=\tilde{R}_2(I,s,\theta),$
then from the definition of $\bar{r}_1$ ,we can know that, $\bar{r}_1$ satisfying
\beq\label{barr1}
c\cdot I^{\frac{1}{p}-i}\leq|{\bar{r}_1}^{(i)}(I)|\leq C \cdot I^{\frac{1}{p}-i},
\eeq
$\bar{r}_2$ have the same estimate with $\tilde{R}_2$ in lemma{\ref{lemma3.2}},i.e.
\beq\label{barr2}
|D_{I}^iD_{s}^j\bar{r}_2)|\le C\cdot
I^{-\frac{i}{q}}, \quad 0\le i+j\le 21. \eeq
 \beq \label{[barr2]}
|D_{I}^i[\bar{r}_2]_0|\le
C\cdot (I^{-i}+I^{-\delta_1-\frac{i}{q}}),\quad 0\le i\le 21. \eeq

 \vskip 0.5cm

\noindent The following results are similarity to \cite{[JPW]},here for the convenience  of readers, we still give the proof of these lemmas.
\begin{Lemma}
Let $0<\delta_1<\frac{1}{10}$ be a constant. Consider the Hamiltonian
\beq\label{lemmabarr} \bar{r}(I,s,\theta)=\omega
I+\bar{r}_1(I)+\mathscr{R}(I,s,\theta), \eeq where
 $\mathscr{R}$ satisfies
\beq\label{mathscrR}|D_{I}^iD_{s}^j\mathscr{R}|\le C\cdot
I^{-\varepsilon-\frac{i}{q}} \eeq for $0\le i+j\le l$ with
$\varepsilon\ge 0$.

 Then there exists a canonical transformation
$\Phi_3$ of the form:
$$\Phi_3:\quad \left\{\begin{array}{ll}
I&=\varrho+u_3(\varrho, \varsigma,\t)\\
s&=\varsigma+v_3(\varrho, \varsigma,\t)
\end{array}\right.
$$
such that
 the
system with Hamiltonian (\ref{lemmabarr}) is transformed into the following one \beq
\label{hatr} \hat{r}(\varrho,\varsigma,\theta)=\omega
\varrho+\hat{r}_1(\varrho)+\mathscr{R}_1(\varrho,\varsigma,\theta),
\eeq where
$\hat{r}_1(\varrho)=\bar{r}_1(\varrho)+[\mathscr{R}]_0(\varrho)$ with $[\mathscr{R}]_0(\varrho)=(\frac{1}{2\pi_p})^2\int_0^{2\pi_p}\!\!\int_0^{2\pi_p} \mathscr{R}(\varrho, \tau,\t)d\tau d\t$ and
$\mathscr{R}_1$ satisfies
\beq\label{mathscrR1define}|D_{\varrho}^iD_{\varsigma}^j\mathscr{R}_1|\le
C\cdot \varrho^{-\varepsilon-\frac{1}{q}-\frac{i}{q}},\quad 0\le  i+j\le
l-3.\eeq
\end{Lemma}

\Proof We will prove this lemma by means of Principle Integral method instead of Fourier series method.
Let $\Phi_3$ be of the following form:
$$
I=\varrho+\frac{\partial S_3}{\partial \tau}(\varrho,s,\theta),\quad \varsigma=s+\frac{\pa S_3}{\pa \varrho}(\varrho,s,\theta),
$$
where the generating function $S_3(\varrho,s,\theta)$ satisfies $S_3(\varrho,s+2\pi_p,\theta)=S_3(\varrho,s,\theta+2\pi_p)=S_3(\varrho,s,\theta)$ and will be determined later.

Then the transformed Hamiltonian is
$$\begin{array}{ll}
\hat{r}&=\omega (\varrho+\frac{\partial S_3}{\partial s})+\bar{r}_1(\varrho+\frac{\partial S_3}{\partial s})+\mathscr{R}(\varrho+\frac{\partial S_3}{\partial s},s,\theta)
+\frac{\partial S_3}{\partial \theta}\\
\\
&=\omega \varrho+\bar{r}_1(\varrho)+[\mathscr{R}]_0(\varrho)+\omega\frac{\partial S_3}{\partial s}+\frac{\partial S_3}{\partial \t}+R+\mathscr{R}_1,\end{array}$$
where
$$R=\mathscr{R}(\varrho,s,\theta)-[\mathscr{R}]_0(\varrho)$$ and
\beq\label{mathscrR1}
\mathscr{R}_1=\int_0^1\bar{r}'_1(\varrho+\lambda \frac{\partial
S_3}{\partial s})\frac{\partial S_3}{\partial s}d\l
+\int_0^1\frac{\pa \mathscr{R}}{\pa I}(\varrho+\lambda
\frac{\partial S_3}{\partial
s},s,\theta)\frac{\partial S_3}{\partial
s}d\l. \eeq
Obviously, it holds that \beq \label{R}
(\frac{1}{2\pi_p})^2\int_0^{2\pi_p}\int_0^{2\pi_p}R(\varrho,s,\theta)ds
d\t=0. \eeq

 Now we
determine the periodic function $S_3$ by the following equation
\beq\label{S3} \omega \frac{\partial S_3}{\partial
s}(\varrho,s,\theta)+\frac{\partial S_3}{\partial
\t}(\varrho,s,\theta)+R(\varrho,s,\theta)=0, \eeq
whose characteristic equation is
$$
\frac{ds}{\omega}=\frac{d\t}{1}=\frac{dS_3}{-R(\varrho,s,\theta)}.
$$
Obviously, the characteristic equation possesses two independent Principle Integrals as follows:
$$s-\o\t=c_1$$
and
$$
S_3+\int_0^{\t}R(\varrho,s-\omega \t+\omega \phi,\phi)d\phi=c_2.
$$
Thus the solution of (\ref{S3}) is of the form:
\beq\label{S_3}S_3(\varrho, s,\t)=-\int_0^{\t}R(\varrho,s-\omega \t+\omega\phi,\phi)d\phi+\Omega(\varrho, s-\omega\t)\eeq
with $\Omega$ a differentiable function determined later.

To ensure $S_3$ be $2\pi_p$-periodic on $s$ and $\t$, $\Omega$ must be $2\pi_p$-periodic on the second variable, that is
$\Omega(\varrho, x+2\pi_p)=\Omega(\varrho, x)$.
Then by direct computation, we obtain that $S_3$ is $2\omega\pi_p$-periodic on $s$.

 Next we determine $\Omega$ by the periodicity of $S_3$ on $\t$.

Let $J(\varrho,x)=-\int_0^{2\pi_p}R(\varrho,x+\omega\phi,\phi)d\phi$. Then we have
$$
\begin{array}{ll}
&S_3(\varrho, s,\t+2\pi_p)=-\int_0^{\t+2\pi_p}R(\varrho,s-\omega (\t+2\pi_p-\phi),\phi)d\phi+\Omega(\varrho, s-\omega (\t+2\pi_p))\\
&=J(\varrho,s-\omega (\t+2\pi_P))-\int_{2\pi_p}^{2\pi_p+\t}R(\varrho,s-\omega (\t+2\pi_p-\phi),\phi)d\phi+\Omega(\varrho, s-\omega (\t+2\pi_p)).
\end{array}
$$
On the other hand, from $R(\varrho,s,\phi+2\pi_p)=R(\varrho,s,\phi)$ we have
$$
\int_{2\pi_p}^{2\pi_p+\t}R(\varrho,s-\omega (\t+2\pi_p-\phi),\phi)d\phi=\int_0^{\t}R(\varrho,s-\omega (\t-\phi),\phi)d\phi,
$$
which implies that
\beq\label{S_32pi}
S_3(\varrho, s,\t+2\pi_p)=J(\varrho,s-\omega (\t+2\pi_p))-\int_0^{\t}R(\varrho,s-\omega (\t-\phi),\phi)d\phi+\Omega(\varrho, s-\omega (\t+2\pi_p)).
\eeq
Setting $S_3(\varrho, s,\t+2\pi_p)=S_3(\varrho, s,\t)$, it follows from (\ref{S_3}) and (\ref{S_32pi}) that
$$
J(\varrho, s-\omega (\t+2\pi_p))+\Omega(\varrho, s-\omega(\t+2\pi_p))-\Omega(\varrho, s-\omega\t)=0,
$$
or equivalently,
\beq\label{Homologous}
J(\varrho, x)=\Omega(\varrho, x+x_0)-\Omega(\varrho, x),
\eeq
where $x=s-\omega (\t+2\pi_p)$ and $x_0=2\omega\pi_p$.

From (\ref{R}) and the definition of $J$, we have
$$
\int_0^{2\pi_p}J(\varrho,x)dx=-\int_0^{2\pi_p}\int_0^{2\pi_p}R(\varrho,x+\o\phi,\phi)dx d\phi
=-\int_0^
{2\omega\pi_p}\int_0^{2\pi_p}R(\varrho,x,\phi)dx d\phi=0.$$
Thus we assume $J(\varrho, x)=\sum_{0\not=k\in Z} J_k(\varrho) e^{i\lambda k x}$ and $\Omega(\varrho, x)=\sum_{0\not=k\in Z} \Omega_k(\varrho) e^{i\lambda k x}$,where $\lambda=\pi/\pi_p$.
Then the homological equation (\ref{Homologous}) implies that
$$
\Omega_k=\frac{J_k}{e^{i\lambda k x_0}-1},\quad k\not=0.
$$
The definition of $J(\varrho,x)$ implies that $J(\varrho,x)$ is $C^l$ on $x$. Thus it holds that
\beq\label{Jk}
|J_k|\le C\cdot\|J(\cdot,x)\|_{C^l}\cdot |k|^{-l},\quad k\not=0.
\eeq
From the Diophantine condition (\ref{Diophantine}), we have that
\beq\label{smalldivisor}
|e^{i\lambda k x_0}-1|\ge 2\pi\gamma|k|^{-\tau},\quad k\not=0.
\eeq
Combining (\ref{Jk}) and (\ref{smalldivisor}), we obtain that
$$
|\Omega_k|\le C\cdot \|J(\cdot,x)\|_{C^l}\cdot|k|^{\tau-l},\quad k\not=0,
$$
which implies $\Omega$ is well-defined and $C^{l-3}$ on $x$ since $1<\tau<2$.

For the definition of $\O$ and (\ref{mathscrR}), we have that
$$|D_{\varrho}^iD_{x}^j{\Omega}|\le C\cdot
\varrho^{-\varepsilon-\frac{i}{q}},\quad 0\le i+j\le l-2,$$
which together with (\ref{mathscrR}) and (\ref{S_3}) implies
\beq\label{S_3inequality}|D_{\varrho}^iD_{\tau}^j{S_3}|\le C\cdot
\varrho^{-\varepsilon-\frac{i}{q}},\quad 0\le i+j\le l-2.\eeq
Thus we obtain (\ref{mathscrR1define}) from (\ref{mathscrR1}) and (\ref{S_3inequality}) and the proof is completed.
\qedbox


By lemma 3.2 and the repeated use of lemma 3.3, we have the
following result.
\begin{Corollary}
There exists a canonical transformation $\Phi_4$ of the form:
$$\Phi_4:\quad \left\{\begin{array}{ll}
I&=\zeta+u_4(\zeta, \eta,\t)\\
s&=\eta+v_4(\zeta, \eta,\t)
\end{array}\right.
$$
such that the
system with Hamiltonian (\ref{barr}) is transformed into the following one \beq
\label{mathfrakr} \mathfrak{r}(\zeta,\eta,\theta)=\omega
\zeta+\mathfrak{r}_1(\zeta)+\mathfrak{r}_2(\zeta,\eta,\theta), \eeq
where $\mathfrak{r}_1=\bar{r}_1+[\bar{r}_2]_0$ with $\bar{r}_1$,
$[\bar{r}_2]_0$ satisfying (\ref{barr1}), (\ref{[barr2]}), and
$\mathfrak{r}_2$ satisfies
\beq\label{mathfrakr2}|D_{\zeta}^iD_{\eta}^j\mathfrak{r}_2|\le
C\cdot \zeta^{-2-\frac{i}{q}} \eeq for $0\le i+j\le 5.$

\end{Corollary}
\section{Proof of theorem 1}
In order to apply Moser's small twist theorem, we need to calculate the pontcare\'e mapping of the Hamiltonian system with the Hamiltonian (\ref{mathfrakr}). So in this section, we first give the expression of the Poincar\'e mapping. And then we will use Moser's small twist theorem to prove Theorem 1.

From corollary 3.1, it follows that the Hamiltonian system with the Hamiltonian (\ref{mathfrakr}) is of
the form: \beq\label{latestsystem} \left\{\begin{array}{ll}
\frac{d\eta}{d\t}&=\o+\mathfrak{r}'_1(\zeta)+\frac{\pa
\mathfrak{r}_2}{\pa \zeta}(\zeta,\eta,\theta)\\
\frac{d\zeta}{d\theta}&=-\frac{\pa \mathfrak{r}_2}{\pa
\eta}(\zeta,\eta,\theta),
\end{array}\right.
\eeq where
$\mathfrak{r}_1(\zeta)=\bar{r}_1(\zeta)+[\bar{r}_2]_0(\zeta)$
satisfying (\ref{barr1}) and (\ref{[barr2]}),
$\mathfrak{r}_2(\zeta,\eta,\theta)$ satisfies (\ref{mathfrakr2}).


Thus the Poincar\'e map of the equation (\ref{latestsystem}) is of
the form: \beq\label{Poincare} P:\ \left\{\begin{array}{ll}
\eta(2\pi_p)&=2\pi_p\o+\eta+\alpha(\zeta)+F_1(\zeta,\eta),\\
\zeta(2\pi_p)&=\zeta+F_2(\zeta,\eta).
\end{array}\right.
\eeq
where $F_1(\zeta,\eta)=\int_0^{2\pi_p} {\frac{\pa
\mathfrak{r}_2}{\pa \zeta}(\zeta,\eta,\theta)}d\theta$, $F_2(\zeta,\eta)=-\int_0^{2\pi_p} {\frac{\pa
\mathfrak{r}_2}{\pa \eta}(\zeta,\eta,\theta)}d\theta$,$\alpha(\zeta)=\mathfrak{r}'_1(\zeta)$,
and from the definition of $\mathfrak{r}_1$, (\ref{barr1}), (\ref{[barr2]}) and (\ref{mathfrakr2}),  we have that \beq\label{alpha}
\alpha(\zeta)=\alpha_1(\zeta)+\alpha_2(\zeta) \eeq with
\beq\label{alphadetail}\begin{array}{ll} &|\alpha^{(i)}_1(\zeta)|\ge
c\cdot \zeta^{-\frac{1}{q}-i},\\
 |\alpha^{(i)}_1(\zeta)|\le
&C\cdot \zeta^{-\frac{1}{q}-i},\quad |\alpha^{(i)}_2(\zeta)|\le
C\cdot \zeta^{-\delta_1-\frac{1}{q}-\frac{i}{q}},\quad 0\le i\le
4\end{array} \eeq \noindent and \beq\label{F1F2}
|D_{\zeta}^iD_{\eta}^jF_k(\zeta,\eta)|\le C\cdot \zeta^{-2-{\frac{i}{q}}},\quad
0\le i+j\le 4,\ k=1, 2,\eeq
where $\alpha_1(\zeta)=\bar{r}'_1(\zeta)$,$\alpha_2(\zeta)={[\bar{r}_2]}'_0(\zeta)$.

 According to (\ref{alphadetail}), we can know that the following case is possible, that is, the function $\alpha(\zeta)$ may be not monotone. In order to find
 a monotone interval for $\alpha(\zeta)$, we consider the interval
 $[2\zeta_0, 3\zeta_0]$ with $\zeta_0\gg 1$. By (\ref{alpha}) and
 (\ref{alphadetail}), we have that the set $\alpha([\frac{9}{4}\zeta_0, \frac{11}{4}\zeta_0])$
 covers some interval with length longer than $c\cdot
 \zeta_0^{-\frac{1}{q}}$. Therefor by Mean Value theorem of Differentials,
 there exists some point $\zeta^*\in[\frac{9}{4}\zeta_0,
 \frac{11}{4}\zeta_0]$ such that $|\alpha'(\zeta^*)|\ge c\cdot
 \zeta_0^{-\frac{1+q}{q}}$.

 What's more, (\ref{alphadetail})
 implies
 $|\alpha''(\zeta)|\le C\cdot \zeta^{-\frac{1+q}{q}-\delta_1}$. Thus for
 each $\zeta\in [\zeta^*,
 \zeta^*+{\zeta_0}^{\frac{\delta_1}{q}}]$, we have
 \beq\label{alpha'}
|\alpha'(\zeta)|\ge c\cdot
 {\zeta_0}^{-\frac{1+q}{q}}.
 \eeq

In the next, we give the following scale transformation :
 \beq\label{scale}
\alpha(\zeta)-\alpha(\zeta^*)=\zeta_0^{-\frac{1+q}{q}}\nu, \quad
\nu\in [2, 3].
 \eeq
Then we have the following Poincar\'e mapping:
\beq\label{tildep} \tilde{P}:\quad\left\{\begin{array}{ll}
\eta(2\pi_p)&=2\pi_p\o+\alpha(\zeta^*)+\eta+\zeta_0^{-\frac{1+q}{q}}\nu+\tilde{F}_1(\nu,\eta),\\
\nu(2\pi_p)&=\nu+\tilde{F}_2(\nu,\eta),
\end{array}\right.
\eeq where \beq\label{tildeF1F2}
 \tilde{F}_1(\nu,\eta)=F_1(\zeta(\nu), \eta), \quad
\tilde{F}_2(\nu,\eta)=\zeta_0^{\frac{1+q}{q}}(\alpha(\zeta(\nu)+F_2(\zeta(\nu),\eta))-\alpha(\zeta(\nu)))
 \eeq with $\zeta(\nu)$ determined by (\ref{scale}).

From (\ref{alphadetail}), (\ref{alpha'}) and (\ref{scale}), we see
that \beq\label{derivative} |\zeta^{(i)}(\nu)|\le C,\quad 1\leq i\le 4,
\eeq which together with (\ref{F1F2}) and (\ref{tildeF1F2}) implies
\beq\label{tildeF1F21} |D_{\nu}^iD_{\eta}^j\tilde{F}_1|\le C\cdot
\zeta_0^{-2},\quad |D_{\nu}^iD_{\eta}^j\tilde{F}_2|\le C\cdot
\zeta_0^{-2}, \quad 0\le i+j\le 4. \eeq

 What's more, the mapping  $\tilde{P}$ of the Hamiltonian system (\ref{mathfrakr}) is
time $2\pi_p$ mapping , so it is area-preserving. And further it possesses the intersection property in the
annulus $[2,3]\times \mathbb{S}_p$, this is to say, if $\Gamma$ is an embedded
circle in $[2,3]\times \mathbb{S}_p$ homotopic to a circle $\nu =$ constant,
then $\tilde{P}(\Gamma)\cap\Gamma\ne\emptyset$. The proof can be
found in \cite{[DZ]}.

For the mapping $\tilde{P}$, all the conditions of Moser's small twist theorem \cite{[Mos]} have been verified.
Consequently, if $\zeta_0\gg 1$, then there exists an invariant curve $\Gamma$ of $\tilde{P}$
surrounding $\nu\equiv 1$ . This implies that the Poincar\'e mapping of the system
(\ref{mathfrakr}) indeed processes invariant curves. Retracting the sequence of transformations back to the original system, we conclude that there exist invariant curves of the Poincar\'e mapping of the original system (\ref{tag(1.10)}). And those curves surround the origin $(x,y)=(0,0)$ and at the same time are
arbitrarily far from it. This completes the proof of Theorem 1.


\end{document}